\documentclass[11pt]{article}%
\usepackage{amsmath}
\usepackage{amssymb}
\usepackage{latexsym}
\usepackage{amsfonts}
\usepackage{graphicx}
\usepackage{color}
\usepackage{tikz}
\usepackage{lineno}%
\newtheorem{prelem}{{\bf Theorem}}

\newtheorem{thm}{Theorem}
\newtheorem{ob}[thm]{Observation}
\newtheorem{lemm}[thm]{Lemma}

\newtheorem{cor}[thm]{Corollary}
\newtheorem{pro}[thm]{Proposition}
\newtheorem{prop}[thm]{Proposition}
\newtheorem{obs}[thm]{Observation}

\DeclareMathOperator{\opt}{\emph{Opt}}

\begin{document}

\title{Total Roman $\{2\}$-Dominating functions in Graphs}
\date{}
\author{$^{(1)}$H. Abdollahzadeh Ahangar, $^{(2)}$M. Chellali,\\$^{(3)}$S.M. Sheikholeslami and $^{(4)}$J.C. Valenzuela-Tripodoro\vspace
{7.5mm} \\$^{(1)}${\small Department of Mathematics}\\{\small Babol Noshirvani University of Technology}\\{\small Shariati Ave., Babol, I.R. Iran, Postal Code: 47148-71167}\\{\small ha.ahangar@nit.ac.ir\vspace{2mm}}\\$^{(3)}${\small LAMDA-RO Laboratory, Department of Mathematics}\\{\small University of Blida}\\{\small Blida, Algeria}\\{\small m\_chellali@yahoo.com\vspace{2mm}}\\$^{(4)}${\small Department of Mathematics }\\{\small Azarbaijan Shahid Madani University}\\{\small Tabriz, Iran}\\{\small s.m.sheikholeslami@azaruniv.ac.ir\vspace{2mm} } \\$^{(5)}${\small Department of Mathematics, University of C\'adiz, Spain.} \\{\small jcarlos.valenzuela@uca.es\vspace{5mm} }}
\maketitle
\begin{abstract}
A Roman $\{2\}$-dominating function (R2F) is a function $f:V\rightarrow
\{0,1,2\}$ with the property that for every vertex $v\in V$ with $f(v)=0$
there is a neighbor $u$ of $v$ with $f(u)=2$, or there are two neighbors $x,y$
of $v$ with $f(x)=f(y)=1$. A total Roman $\{2\}$-dominating function (TR2DF)
is an R2F $f$ such that the set of vertices with $f(v)>0$ induce a subgraph
with no isolated vertices. The weight of a TR2DF is the sum of its function
values over all vertices, and the minimum weight of a TR2DF of $G$ is the
total Roman $\{2\}$-domination number $\gamma_{tR2}(G).$ In this paper, we
initiate the study of total Roman $\{2\}$-dominating functions, where
properties are established. Moreover, we present various bounds on the total
Roman $\{2\}$-domination number. We also show that the decision problem
associated with $\gamma_{tR2}(G)$ is NP-complete for bipartite and chordal graphs.
{Moreover, we show that it is possible to
compute this parameter in linear time for bounded clique-width graphs (including trees).}


\end{abstract}

\section{Introduction}

In this paper, $G$ is a simple graph with vertex set $V=V(G)$ and edge set
$E=E(G)$. The order $|V|$ of $G$ is denoted by $n(G)$. For every vertex $v\in
V$, the \emph{open neighborhood} $N(v)$ is the set $\{u\in V(G):uv\in E(G)\}$
and the \emph{closed neighborhood} of $v$ is the set $N[v]=N(v)\cup\{v\}$. The
\emph{degree} of a vertex $v\in V$ is $\deg_{G}(v)=|N(v)|$. The \emph{minimum}
and \emph{maximum degree} of a graph $G$ are denoted by $\delta=\delta(G)$ and
$\Delta=\Delta(G)$, respectively. A \textit{leaf} of $G$ is a vertex of degree
one, while a \textit{support vertex} of $G$ is a vertex adjacent to a leaf. A
support vertex is said to be \textit{weak }(resp\textit{. strong})\textit{ }if
it is adjacent to exactly one leaf (resp. \textit{at least two leaves}).

Let $P_{n}$, $C_{n}$ and $K_{n}$ be the \emph{path}, \emph{cycle} and
\emph{complete graph} of order $n$ and $K_{p,q}$ the complete bipartite graph
with one partite set of cardinality $p$ and the other of cardinality $q$. The
\emph{complement} of a graph $G$ is denoted by $\overline{G}$. The
\textit{join} of two graphs $G$ and $H$, denoted by $G\vee H$, is a graph
obtained from $G$ and $H$ by joining each vertex of $G$ to all vertices of
$H.$ A \emph{tree} is an acyclic connected graph. A \emph{double star} is a
tree containing exactly two vertices that are not leaves. A double star with
respectively $p$ and $q$ leaves attached at each support vertex is denoted by
$S_{p,q}.$ {The \textit{corona} of a graph $H$
is the graph obtained from $H$ by appending a
vertex of degree 1 to each vertex of $H$. The \textit{distance}
$d_{G}(u,v)$ between two vertices $u$ and $v$ in a connected
graph $G$ is the length of a shortest $u-v$ path in $G.$ The
\emph{diameter} of $G$, denoted by $\mathrm{diam}(G)$, is the
maximum value among distances between all pair of vertices of $G$.}

A subset $S\subseteq V$ is a \emph{dominating set} if every vertex in
{$V-S$} has a neighbor in $S$, and $S$ is a \emph{total dominating
set, }abbreviated TSD, if every vertex in $V$ has a neighbor in $S$, that is,
$N(v)\cap S\neq\emptyset$ for all $v\in V$. The \textit{domination number}
$\gamma(G)$ of a graph $G$ is the minimum cardinality of a dominating set of
$G$, and the \emph{total domination number} $\gamma_{t}(G)$ is the minimum
cardinality of a TDS of $G.$

For a graph $G$ and a positive integer $k$, let $f:V(G)\rightarrow
\{0,1,2,...,k\}$ be a function, and let $(V_{0},V_{1},V_{2},\ldots,V_{k})$ be
the ordered partition of $V=V(G)$ induced by $f$, where $V_{i}=\{v\in
V:f(v)=i\}$ for $i\in\{0,1,\ldots,k\}$. There is a $1$-$1$ correspondence
between the functions $f:V\rightarrow\{0,1,2,...,k\}$ and the ordered
partitions $(V_{0},V_{1},V_{2},\ldots,V_{k})$ of $V$, so we will write
$f=(V_{0},V_{1},V_{2},\ldots,V_{k})$. {The weight of }$f$ is the value
$f(V(G))=\sum_{u\in V(G)}f(u).$

A function $f:V(G)\rightarrow\{0,1,2\}$ is a \textit{Roman dominating
function, }abbreviated RDF, on $G$ if every vertex $u\in V$ for which $f(u)=0$
is adjacent to at least one vertex $v$ for which $f(v)=2$. The \textit{Roman
domination number} $\gamma_{R}(G)$ is the minimum weight of a RDF on $G$.
Roman domination was introduced by Cockayne et al. in \cite{CDHH} and was
inspired by the work of ReVelle and Rosing~\cite{ReVelle}, and
Stewart~\cite{Stewart}.

{ The definition of Roman dominating functions was motivated by an
article in {\it Scientific American} by Ian Stewart entitled ``Defend
the Roman Empire"~\cite{Stewart} and suggested even earlier by
ReVelle~\cite{Re97}. Each vertex in our graph represents a location
in the Roman Empire. A location (vertex $v$) is considered {\em
unsecured} if no legions are stationed there (i.e., $f(v) = 0$) and
{\em secured} otherwise (i.e., if $f(v) \in \{1,2\})$. An unsecured
location (vertex $v$) can be secured by sending a legion to $v$ from
a neighboring location (an adjacent vertex $u$). But Constantine the
Great (Emperor of Rome) issued a decree in the 4th century A.D. for
the defense of the regions, where a legion cannot be sent from a secured
location having only one legion stationed there to an unsecured location,
for otherwise it leaves that location unsecured. Thus, two legions must be stationed at a
location ($f(v) = 2$) before one of the legions can be sent to a neighboring location. In this way, Emperor Constantine the Great can
defend the Roman Empire. Since it is expensive to maintain a legion
at a location, the Emperor would like to station as few legions as
possible, while still defending the Roman Empire. A Roman dominating
function of weight $\gamma_R(G)$ corresponds to such an optimal assignment
of legions to locations.}

In \cite{chha}, Chellali et al. introduced the Roman $\{2\}$-domination
{(called in \cite{KM} and elsewhere Italian domination) }defined as follows: a
\textit{Roman }$\{2\}$\textit{-dominating function} is a function
$f=(V_{0},V_{1},V_{2})$ with the property that for every vertex $v\in V_{0}$
there is a vertex $u\in N(v)$, with $u\in V_{2}$ or there are two vertices
$x,y\in N(v)$ with $x,y\in V_{1}$. The \textit{Roman }$\{2\}$%
\textit{-domination number }$\gamma_{\{R2\}}(G)$ is the minimum weight of a
Roman $\{2\}$-dominating function on $G.$

{
There are many papers in the literature devoted to the study of Roman domination
type problems and its variations. One of the questions that arise naturally
when a Roman domination type problem is studied is to focus on its complexity and
algorithmic aspects. In 2008, Liedloff et all \cite{LKLP} showed,
among other results, that it is possible to compute the Roman domination number
of a graph with bounded cliquewidth in linear time. Clearly, this implies that there
exists algorithms for computing the Roman domination number of trees in linear time.
Chellali et al. \cite{chha} proved that Roman $\{2\}$-domination problem is NP-complete for bipartite graphs while Chen and Lu \cite{CL} recently showed it is NP-complete even when restricted to split graphs. Moreover, the authors \cite{CL} presented a linear time algorithm to obtain the Roman $\{2\}$-domination number of a block graph.
}

In this paper, we initiate the study of the total version of Roman
$\{2\}$-dominating function.\ A \textit{total Roman }$\{2\}$%
\textit{-dominating function}, abbreviated TR2DF, is a Roman $\{2\}$%
-dominating function $f=(V_{0},V_{1},V_{2})$ such that the subgraph induced by
$V_{1}\cup V_{2}$ has no isolated vertices. The \textit{total Roman }%
$\{2\}$\textit{-domination number }$\gamma_{tR2}(G)$ is the minimum weight of
a TR2DF on $G.$ {A TR2DF on }$G$ with weight $\gamma_{tR2}(G)$ is called an
$\gamma_{tR2}(G)$-function. {Total Roman }$\{2\}$-domination number is
well-defined for all graphs $G$ with no isolated vertices since assigning a
$1$ to every vertex of $G$ provides a TR2DF of $G$. {Hence for all graphs }$G$
of order $n$ with $\delta(G)\geq1,$ $2\leq\gamma_{tR2}(G)\leq n.$
{ We present various bounds on the total Roman $\{2\}$domination number and several
properties are established. We show that the decision problem associated with $\gamma_{tR2}(G)$ is
NP-complete for bipartite and chordal graphs. Moreover, we show that it is possible to
compute this parameter in linear time for bounded clique-width graphs (including trees).
}

We note that throughout this paper, we only consider nontrivial connected graphs that
{we will call \textit{ntc graphs}.}

\section{Preliminary results}

{We begin by giving some properties of total Roman }$\{2\}$-dominating
functions. The following two facts lead to our first observation. Clearly
assigning a $2$ to every vertex in a minimum total dominating set of a ntc
graph $G$ and a $0$ to the remaining vertices of $G$ provides a TR2DF. Also,
for every TR2DF $f=(V_{0},V_{1},V_{2})$ the set $V_{1}\cup V_{2}$ total
dominates $V(G).$

\begin{obs}
\label{lower and upper}For every ntc graph $G,$
\[
\gamma_{t}(G)\leq\gamma_{tR2}(G)\leq2\gamma_{t}(G).
\]

\end{obs}

{Note that the lower bound of }Observation \ref{lower and upper} is attained
for $K_{2}\vee\overline{K_{n-2}}$ {while the upper bound is attained for the
double star }$S_{3,3}.$

{It is well-known that }$\gamma_{t}(G)\leq2\gamma(G)$ {for every ntc graph
}$G.$ {So by Observation \ref{lower and upper}, we will have }$\gamma
_{tR2}(G)\leq4\gamma(G).$ {Our next result improves this upper bound to
}$\gamma_{tR2}(G)\leq3\gamma(G).$ We need the following result due to
Bollob\'{a}s and Cockayne \cite{bc}. {If $S$ is a set
of vertices, then we say that a vertex $v$ is a \emph{private neighbor} of a
vertex $u\in S$ (with respect to $S$) if $N[v]\cap S=\{u\}$. The
\textit{external private neighborhood} \textrm{epn(}$u,S)$ of $u$ {with
respect to }${S}${ consists of those private neighbors of }${u}${ in }$V-S.$
{For a TR2DF $f=(V_{0},V_{1},V_{2})$ of an ntc graph, let {$V_{02}=\{w\in
V_{0}:N(w)\cap V_{2}\neq\emptyset\}$ and $V_{01}=V_{0}-V_{02}.$}}}

\begin{thm}
[Bollob\'{a}s and Cockayne \cite{bc}]\label{bc-minimality}If $G$ is a graph
without isolated vertices, then $G$ has a minimum dominating set $D$ such that
for all $d\in D$, there exists a neighbor {$f(d)\in V-D$} of $d$
such that $f(d)$ is not a neighbor of any vertex {$x\in D-\{d\}.$}
\end{thm}

\begin{prop}
\label{3gamma}For every ntc graph $G,$ $\gamma_{tR2}(G)\leq3\gamma(G).$
\end{prop}

\textbf{Proof. }Let $D$ be a minimum dominating set of $G$ satisfying the
property of Theorem~\ref{bc-minimality}. Since each vertex of $D$ has an
external private {neighbor in $V-D,$} let $W$ be a subset of
{$V-D$} formed by the private neighbors chosen so that each vertex
of $D$ has exactly one external private neighbor in $D.$ Clearly $\left\vert
W\right\vert =\left\vert D\right\vert .$ Now define the function $f$ as
follows: $f(x)=2$ for all $x\in D,$ $f(x)=1$ for all $x\in W$, and $f(x)=0$
otherwise. Clearly $f$ is a TR2DF of $G$ of weight $2\left\vert D\right\vert
+\left\vert W\right\vert =3\gamma(G),$ and thus $\gamma_{tR2}(G)\leq
3\gamma(G).$ $\ \Box$

{For the sharpness of the bound in Proposition \ref{3gamma},
consider the tree $T$ obtained from two stars $K_{1,4}$ by
adding an edge between a leaf of one star to a leaf of the other star. The
next observation is straightforward and is tight for double stars. }

\begin{obs}
For every ntc graph $G,$ $\gamma_{R2}(G)\leq\gamma_{tR2}(G)$.
\end{obs}

\begin{prop}
\label{po}Let $G$ be an ntc graph. Then for every $\gamma_{tR2}(G)$-function
$f=(V_{0},V_{1},V_{2})$ such that {$|V_{2}|$ is as small as possible, we have:
}

\begin{description}
\item[i) ] Each vertex in $V_{2}$ has at least two private neighbors in
$V_{0}$ with respect to $V_{2}$.

\item[ii) ] $2|V_{2}|\leq|V_{02}|.$
\end{description}
\end{prop}

\textbf{Proof. } i) Suppose there exists a vertex $v\in V_{2}$ with
{at most one private neighbor} in $V_{0}$ with respect to $V_{2}$.
Then reassigning $v$ and its private neighbor (if any) the value $1$ instead
of $2$ and $0,$ respectively provides a $\gamma_{tR2}(G)$-function with less
vertices assigned a $2,$ which contradicts the choice of $f.$

ii)- follows from (i). $\ \Box$

\begin{prop}
\label{emptyV2}Let $G$ be an ntc graph with maximum degree $\Delta\leq2$.
{Then there exists} a $\gamma_{tR2}(G)$-function $f=(V_{0}%
,V_{1},V_{2})$ such that $V_{2}=\emptyset.$
\end{prop}

\textbf{Proof. }{Among all }$\gamma_{tR2}(G)$-functions, let $f=(V_{0}%
,V_{1},V_{2})$ be a one such that {$|V_{2}|$ is as small as possible. If
}$V_{2}\neq\emptyset,$ then by Proposition \ref{po}, every vertex of $V_{2}$
has at least two private neighbors in $V_{0}$ with respect to $V_{2}.$ But
then {since $\Delta\leq2,$ each vertex in }$V_{2}$
would be isolated in $G[V_{1}\cup V_{2}],$ a contradiction. Hence
{$V_{2}=\emptyset$.}
$\hfill\Box$

{Recall that a subset }$S$ of $V$ is a \textit{double dominating set} of\ $G$
if for every vertex $v\in V,$ $\left\vert N[v]\cap S\right\vert \geq2,$ that
is, $v$ is in $S$ and has at least one neighbor in $S$ or $v$ is in $V-S$ and
has at least two neighbors in $S.$ The \textit{double domination number}
$\gamma_{\times2}(G)$ is the minimum cardinality of a double dominating set of
$G$. Double domination was introduced by Harary and Haynes \cite{hv1} who
proved that $\gamma_{\times2}(P_{n})=\lceil\frac{2n+2}{3}\rceil$ and
$\gamma_{\times2}(C_{n})=\lceil\frac{2n}{3}\rceil$. {The following result
shows that the equality between }$\gamma_{tR2}(G)$ and $\gamma_{\times2}(G)$
{occurs }under certain conditions.

\begin{prop}
\label{double}Let $G$ be an ntc graph. If $G$ has a $\gamma_{tR2}(G)$-function
$f=(V_{0},V_{1},V_{2})$ such that $V_{2}=\emptyset$, then $\gamma
_{tR2}(G)=\gamma_{\times2}(G)$.
\end{prop}

\noindent\textbf{Proof. } If $S$ is a double dominating set of $G,$ then
$(V-S,S,\emptyset)$ is clearly a TR2DF on $G,$ and thus $\gamma_{tR2}%
(G)\leq\gamma_{\times2}(G).$ Now if $f=(V_{0},V_{1},V_{2})$ is a $\gamma
_{tR2}(G)$-function such that $V_{2}=\emptyset$, then $V_{1}$ double dominates
$V(G),$ and thus $\gamma_{\times2}(G)\leq\gamma_{tR2}(G).$ Therefore
$\gamma_{tR2}(G)=\gamma_{\times2}(G)$. $\hfill\Box$

The following results are immediate consequences of Propositions
\ref{emptyV2}, \ref{double} and the exact values of the double domination
number of paths and cycles given above.

\begin{prop}
\label{path}For $n\geq2$, $\gamma_{tR2}(P_{n})=\lceil\frac{2n+2}{3}\rceil$.
\end{prop}

\begin{prop}
\label{cycle2}For $n\geq3$, $\gamma_{tR2}(C_{n})=\lceil\frac{2n}{3}\rceil$.
\end{prop}

\section{Complexity}

Our aim in this section is to study the complexity of the following decision
problem, to which we shall refer as TOTAL ROMAN $\{2\}$-DOMINATION:

\bigskip

TOTAL ROMAN $\{2\}$-DOMINATION

\textbf{Instance}: Graph $G=(V,E)$, positive integer $k\leq|V|$.

\textbf{Question}: Does $G$ have a total Roman $\{2\}$-dominating function of
weight at most $k$?

\bigskip

We show that this problem is NP-complete by reducing the well-known
NP-complete problem, Exact-$3$-Cover (X3C), to TOTAL ROMAN $\{2\}$-DOMINATION.

\bigskip

EXACT $3$-COVER (X3C)

\textbf{Instance}: A finite set $X$ with $|X|=3q$ and a collection $C$ of
$3$-element subsets of $X$.

\textbf{Question}: Is there a subcollection $C^{\prime}$ of $C$ such that
every element of $X$ appears in exactly one element of $C^{\prime}$?

\begin{figure}[h]
\centering
\begin{tikzpicture}
\node [draw, shape=circle,fill=white,scale=0.5] (b1) at  (0,0.8) {};
\node [draw, shape=circle,fill=white,scale=0.5] (a2) at  (0,1.5) {};
\node [draw, shape=circle,fill=white,scale=0.5] (a3) at  (-.5,2) {};
\node [draw, shape=circle,fill=white,scale=0.5] (a4) at  (.5,2) {};
\node [draw, shape=circle,fill=white,scale=0.5] (a5) at  (-.8,2.6) {};
\node [draw, shape=circle,fill=white,scale=0.5] (a6) at  (-.5,2.6) {};
\node [draw, shape=circle,fill=white,scale=0.5] (a7) at  (-.2,2.6) {};
\node [draw, shape=circle,fill=white,scale=0.5] (a8) at  (.8,2.6) {};
\node [draw, shape=circle,fill=white,scale=0.5] (a9) at  (.5,2.6) {};
\node [draw, shape=circle,fill=white,scale=0.5] (a10) at  (.2,2.6) {};
\draw(b1)--(a2)--(a3)--(a5)--(a3)--(a6)--(a3)--(a7)--(a3)--(a2)--(a4)--(a8)--(a4)--(a9)--(a4)--(a10);
\node [draw, shape=circle,fill=white,scale=0.5] (b2) at  (2,0.8) {};
\node [draw, shape=circle,fill=white,scale=0.5] (a2) at  (2,1.5) {};
\node [draw, shape=circle,fill=white,scale=0.5] (a3) at  (1.5,2) {};
\node [draw, shape=circle,fill=white,scale=0.5] (a4) at  (2.5,2) {};
\node [draw, shape=circle,fill=white,scale=0.5] (a5) at  (1.2,2.6) {};
\node [draw, shape=circle,fill=white,scale=0.5] (a6) at  (1.5,2.6) {};
\node [draw, shape=circle,fill=white,scale=0.5] (a7) at  (1.8,2.6) {};
\node [draw, shape=circle,fill=white,scale=0.5] (a8) at  (2.8,2.6) {};
\node [draw, shape=circle,fill=white,scale=0.5] (a9) at  (2.5,2.6) {};
\node [draw, shape=circle,fill=white,scale=0.5] (a10) at  (2.2,2.6) {};
\draw(b2)--(a2)--(a3)--(a5)--(a3)--(a6)--(a3)--(a7)--(a3)--(a2)--(a4)--(a8)--(a4)--(a9)--(a4)--(a10);
\node [draw, shape=circle,fill=white,scale=0.5] (b3) at  (4,0.8) {};
\node [draw, shape=circle,fill=white,scale=0.5] (a2) at  (4,1.5) {};
\node [draw, shape=circle,fill=white,scale=0.5] (a3) at  (3.5,2) {};
\node [draw, shape=circle,fill=white,scale=0.5] (a4) at  (4.5,2) {};
\node [draw, shape=circle,fill=white,scale=0.5] (a5) at  (3.2,2.6) {};
\node [draw, shape=circle,fill=white,scale=0.5] (a6) at  (3.5,2.6) {};
\node [draw, shape=circle,fill=white,scale=0.5] (a7) at  (3.8,2.6) {};
\node [draw, shape=circle,fill=white,scale=0.5] (a8) at  (4.8,2.6) {};
\node [draw, shape=circle,fill=white,scale=0.5] (a9) at  (4.5,2.6) {};
\node [draw, shape=circle,fill=white,scale=0.5] (a10) at  (4.2,2.6) {};
\draw(b3)--(a2)--(a3)--(a5)--(a3)--(a6)--(a3)--(a7)--(a3)--(a2)--(a4)--(a8)--(a4)--(a9)--(a4)--(a10);
\node [draw, shape=circle,fill=white,scale=0.5] (b4) at  (6,0.8) {};
\node [draw, shape=circle,fill=white,scale=0.5] (a2) at  (6,1.5) {};
\node [draw, shape=circle,fill=white,scale=0.5] (a3) at  (5.5,2) {};
\node [draw, shape=circle,fill=white,scale=0.5] (a4) at  (6.5,2) {};
\node [draw, shape=circle,fill=white,scale=0.5] (a5) at  (5.2,2.6) {};
\node [draw, shape=circle,fill=white,scale=0.5] (a6) at  (5.5,2.6) {};
\node [draw, shape=circle,fill=white,scale=0.5] (a7) at  (5.8,2.6) {};
\node [draw, shape=circle,fill=white,scale=0.5] (a8) at  (6.8,2.6) {};
\node [draw, shape=circle,fill=white,scale=0.5] (a9) at  (6.5,2.6) {};
\node [draw, shape=circle,fill=white,scale=0.5] (a10) at  (6.2,2.6) {};
\draw(b4)--(a2)--(a3)--(a5)--(a3)--(a6)--(a3)--(a7)--(a3)--(a2)--(a4)--(a8)--(a4)--(a9)--(a4)--(a10);
\node [draw, shape=circle,fill=white,scale=0.5] (b5) at  (8,0.8) {};
\node [draw, shape=circle,fill=white,scale=0.5] (a2) at  (8,1.5) {};
\node [draw, shape=circle,fill=white,scale=0.5] (a3) at  (7.5,2) {};
\node [draw, shape=circle,fill=white,scale=0.5] (a4) at  (8.5,2) {};
\node [draw, shape=circle,fill=white,scale=0.5] (a5) at  (7.2,2.6) {};
\node [draw, shape=circle,fill=white,scale=0.5] (a6) at  (7.5,2.6) {};
\node [draw, shape=circle,fill=white,scale=0.5] (a7) at  (7.8,2.6) {};
\node [draw, shape=circle,fill=white,scale=0.5] (a8) at  (8.8,2.6) {};
\node [draw, shape=circle,fill=white,scale=0.5] (a9) at  (8.5,2.6) {};
\node [draw, shape=circle,fill=white,scale=0.5] (a10) at  (8.2,2.6) {};
\draw(b5)--(a2)--(a3)--(a5)--(a3)--(a6)--(a3)--(a7)--(a3)--(a2)--(a4)--(a8)--(a4)--(a9)--(a4)--(a10);
\node [draw, shape=circle,fill=white,scale=0.5] (b6) at  (10,0.8) {};
\node [draw, shape=circle,fill=white,scale=0.5] (a2) at  (10,1.5) {};
\node [draw, shape=circle,fill=white,scale=0.5] (a3) at  (9.5,2) {};
\node [draw, shape=circle,fill=white,scale=0.5] (a4) at  (10.5,2) {};
\node [draw, shape=circle,fill=white,scale=0.5] (a5) at  (9.2,2.6) {};
\node [draw, shape=circle,fill=white,scale=0.5] (a6) at  (9.5,2.6) {};
\node [draw, shape=circle,fill=white,scale=0.5] (a7) at  (9.8,2.6) {};
\node [draw, shape=circle,fill=white,scale=0.5] (a8) at  (10.8,2.6) {};
\node [draw, shape=circle,fill=white,scale=0.5] (a9) at  (10.5,2.6) {};
\node [draw, shape=circle,fill=white,scale=0.5] (a10) at  (10.2,2.6) {};
\draw(b6)--(a2)--(a3)--(a5)--(a3)--(a6)--(a3)--(a7)--(a3)--(a2)--(a4)--(a8)--(a4)--(a9)--(a4)--(a10);
\node [draw, shape=circle,fill=white,scale=0.5] (c1) at  (.95,-1.2) {};
\node [draw, shape=circle,fill=white,scale=0.5] (c2) at  (3.65,-1.2) {};
\node [draw, shape=circle,fill=white,scale=0.5] (c3) at  (6.35,-1.2) {};
\node [draw, shape=circle,fill=white,scale=0.5] (c4) at  (9.05,-1.2) {};
\draw(c2)--(b1)--(c1)--(b4)--(c4)--(b3)--(c3)--(b2)--(c4);
\draw(c3)--(b6)--(c2)--(b5)--(c1);
\node [draw, shape=circle,fill=white,scale=0.5] (d1) at  (.95,-1.8) {};
\node [draw, shape=circle,fill=white,scale=0.5] (d2) at  (.45,-1.8) {};
\node [draw, shape=circle,fill=white,scale=0.5] (d3) at  (1.45,-1.8) {};
\node [draw, shape=circle,fill=white,scale=0.5] (d4) at  (.95,-2.5) {};
\node [draw, shape=circle,fill=white,scale=0.5] (d5) at  (.45,-3.2) {};
\node [draw, shape=circle,fill=white,scale=0.5] (d6) at  (.95,-3.2) {};
\node [draw, shape=circle,fill=white,scale=0.5] (d7) at  (1.45,-3.2) {};
\draw(d6)--(d4)--(d5)--(d4)--(d7)--(d4)--(d1)--(d3)--(d1)--(d2)--(d1)--(c1);
\node [draw, shape=circle,fill=white,scale=0.5] (d1) at  (3.65,-1.8) {};
\node [draw, shape=circle,fill=white,scale=0.5] (d2) at  (3.15,-1.8) {};
\node [draw, shape=circle,fill=white,scale=0.5] (d3) at  (4.15,-1.8) {};
\node [draw, shape=circle,fill=white,scale=0.5] (d4) at  (3.65,-2.5) {};
\node [draw, shape=circle,fill=white,scale=0.5] (d5) at  (3.15,-3.2) {};
\node [draw, shape=circle,fill=white,scale=0.5] (d6) at  (3.65,-3.2) {};
\node [draw, shape=circle,fill=white,scale=0.5] (d7) at  (4.15,-3.2) {};
\draw(d6)--(d4)--(d5)--(d4)--(d7)--(d4)--(d1)--(d3)--(d1)--(d2)--(d1)--(c2);
\node [draw, shape=circle,fill=white,scale=0.5] (d1) at  (6.35,-1.8) {};
\node [draw, shape=circle,fill=white,scale=0.5] (d2) at  (5.85,-1.8) {};
\node [draw, shape=circle,fill=white,scale=0.5] (d3) at  (6.85,-1.8) {};
\node [draw, shape=circle,fill=white,scale=0.5] (d4) at  (6.35,-2.5) {};
\node [draw, shape=circle,fill=white,scale=0.5] (d5) at  (5.85,-3.2) {};
\node [draw, shape=circle,fill=white,scale=0.5] (d6) at  (6.35,-3.2) {};
\node [draw, shape=circle,fill=white,scale=0.5] (d7) at  (6.85,-3.2) {};
\draw(d6)--(d4)--(d5)--(d4)--(d7)--(d4)--(d1)--(d3)--(d1)--(d2)--(d1)--(c3);
\node [draw, shape=circle,fill=white,scale=0.5] (d1) at  (9.05,-1.8) {};
\node [draw, shape=circle,fill=white,scale=0.5] (d2) at  (8.55,-1.8) {};
\node [draw, shape=circle,fill=white,scale=0.5] (d3) at  (9.55,-1.8) {};
\node [draw, shape=circle,fill=white,scale=0.5] (d4) at  (9.05,-2.5) {};
\node [draw, shape=circle,fill=white,scale=0.5] (d5) at  (8.55,-3.2) {};
\node [draw, shape=circle,fill=white,scale=0.5] (d6) at  (9.05,-3.2) {};
\node [draw, shape=circle,fill=white,scale=0.5] (d7) at  (9.55,-3.2) {};
\draw(d6)--(d4)--(d5)--(d4)--(d7)--(d4)--(d1)--(d3)--(d1)--(d2)--(d1)--(c4);
\end{tikzpicture}
\caption{NP-Completeness for bipartite graphs.}%
\label{fig-1}%
\end{figure}
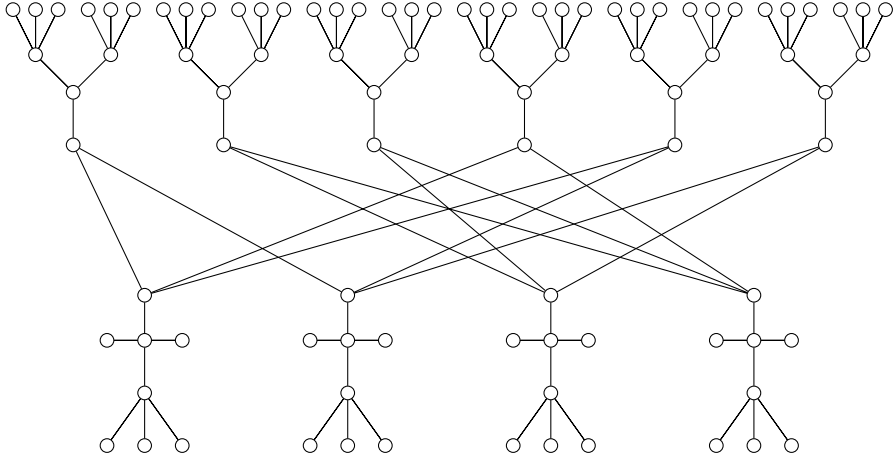

\begin{thm}
\label{complexity_bipartite}TOTAL\ ROMAN $\{2\}$-DOMINATION is NP-Complete for
bipartite graphs.
\end{thm}

\textbf{Proof. }TOTAL\ ROMAN $\{2\}$-DOMINATION is a member of $\mathcal{NP}$,
since we can check in polynomial time that a function $f:V\rightarrow
\{0,1,2\}$ has weight at most $k$ and is a total Roman $\{2\}$-dominating
function. Now let us show how to transform any instance of X3C into an
instance $G$ of TOTAL\ ROMAN $\{2\}$-DOMINATION so that one of them has a
solution if and only if the other one has a solution. Let $X=\{x_{1}%
,x_{2},...,x_{3q}\}$ and $C=\{C_{1},C_{2},...,C_{t}\}$ be an arbitrary
instance of X3C.

For each $x_{i}\in X$, we build a graph $H_{i}$ obtained from a path
$P_{2}:x_{i}$-$y_{i}$ and two stars $K_{1,3}$ centered at $a_{i}$ and $b_{i},$
by adding edges $y_{i}a_{i}$ and $y_{i}b_{i}.$ Hence, each $H_{i}$ has order
$10$. For each $C_{j}\in C$, we build a double star $S_{3,3}$ with support
vertices $u_{j}$ and $v_{j}.$ Let $c_{j}$ be a leaf of the double star
$S_{3,3}.$ Let $Y=\{c_{1},c_{2},...,c_{t}\}$. Now to obtain a graph $G$, we
add edges $c_{j}x_{i}$ if $x_{i}\in C_{j}$. Clearly $G$ is a bipartite graph.
Set $k=4t+16q$. Observe that for every total Roman $\{2\}$-dominating function
$f$ on $G$, each $H_{i}$ has weight at least $5$ and each double star
$S_{3,3}$ has weight at least $4$.

Suppose that the instance $X,C$ of X3C has a solution $C^{\prime}$. We
construct a total Roman $\{2\}$-dominating function $f$ on $G$ of weight $k$.
For each $i,$ assign the value $2$ to $a_{i},b_{i};$ the value $1$ to $y_{i}$
and $0$ to the remaining vertices of $H_{i}.$ For every $j,$ assign the value
$2$ to $u_{j}$ and $v_{j},$ and $0$ to each leaf. In addition, for every
$c_{j}$, assign the value $1$ if $C_{j}\in C^{\prime}$ and the value $0$ if
$C_{j}\notin C^{\prime}.$ Note that since $C^{\prime}$ exists, its cardinality
is precisely $q$, and so the number of $c_{j}$'s with value $1$ is $q$, having
disjoint neighborhoods in $\{x_{1},x_{2},...,x_{3q}\}${. Since
$C^{\prime}$} is a solution for X3C, every vertex in $X$ is adjacent to two
vertices assigned a $1$. Hence, it is straightforward to see that $f$ is a
total Roman $\{2\}$-dominating function with weight $f(V)=4t+q+15q=k.$

Conversely, suppose that $G$ has a total Roman $\{2\}$-dominating function
with weight at most $k$. Among all such functions, let $g=(V_{0},V_{1},V_{2})$
be one such that $\{y_{1},y_{2},...,y_{3q}\}\cap V_{2}$ is as small as
possible. As observed above, since each $H_{i}$ has weight at least $5,$ we
may assume that $g(a_{i})=g(b_{i})=2$ and $g(y_{i})>0$ so that vertices
$a_{i},b_{i}$ are not isolated in the subgraph induced by $V_{1}\cup V_{2}.$
Hence each leaf neighbor of $a_{i}$ or $b_{i}$ is assigned a $0$ under $g.$ If
$g(y_{i})=2$ for some $i,$ then we must have $g(x_{i})=0.$ In this case,
reassigning a $1$ to each of $y_{i}$ and $x_{i}$ instead of $2$ and $0,$
respectively, provides a total Roman $\{2\}$-dominating function $g^{\prime}$
with less vertices $y_{i}$ assigned a $2$ than under $g$, contradicting our
choice of $g.$ Hence $g(y_{i})=1$ for every $i\in\{1,2,...,3q\}.$ On the other
hand, the total weight of all double stars corresponding to
{elements of $C$ is at least $4t$.} In this case, we can assume
that $g(u_{j})=g(v_{j})=2$ and so each leaf neighbor of $u_{j}$ or $b_{j}$ is
assigned a $0$ under $g.$ Note that each $c_{j}$ can be assigned a $0$ since
$g(u_{j})=2.$ Since $w(g)\leq4t+16q$ and the total weight assigned to vertices
of $V(G)-(X\cup Y)$ is $4t+15q,$ we have to assign to vertices of $(X\cup Y)$
weights whose total not exceeding $q$ {in order that} each vertex
$x_{i}\in X$ has either $g(x_{i})>0$ or has two neighbors in $V_{1}.$ Since
$\left\vert X\right\vert =3q,$ it is clear that this is only possible if there
are $q$ vertices of $\{c_{1},c_{2},...,c_{t}\}$ that are assigned a $1.$ Since
each $c_{j}$ has a exactly three neighbors in $\{x_{1},x_{2},...,x_{3q}\},$ we
deduce that $C^{\prime}=\{C_{j}:g(c_{j})=1\}$ is an exact cover for $C$.
$\ \Box$ \medskip

The next result is obtained by using the same proof of Theorem
\ref{complexity_bipartite} on the (same) graph $G$ built for the
transformation by adding all edges between the vertices labelled $c_{j}$
{in order that} the resulting graph is chordal.

\begin{thm}
TOTAL\ ROMAN $\{2\}$-DOMINATION is NP-Complete for chordal graphs.
\end{thm}

{
In the rest of this section, we prove that the decision problem associated to
$\gamma_{tR2}(G)$ can be solved in linear time for the class of graphs with
bounded clique-width, which implies that it also can be computed in linear time
for trees.

We make use of several useful objects and results, which are formally described in \cite{CMR, LKLP},
related to logic structures. Namely, in what follows, a $k$-expression on the vertices of a
graph $G$ with labels $\{ 1,2,\ldots , k \}$ is an expression using the following operations:

\begin{tabular}{rl}
	$\bullet i(x)$ & create a new vertex $x$ with label $i$ \\
	$G_1\oplus G_2$ & create a new graph which is the disjoint union of the graphs $G_i$ \\
	$\eta_{ij}(G)$ & add all edges in $G$ joining vertices with label $i$ with vertices with label $j$ \\
	$\rho_{i\rightarrow j}(G)$ & change the label of all vertices with label $i$ into label $j$
\end{tabular}

\noindent We call the {\it clique-width} of a graph $G$ the minimum positive integer $k$ that is
needed to describe $G$ by means of a $k$-expression. For example, the complete graph $K_3$ with
set of vertices $\{a,b,c\}$ can be described by the following $2$-expression.
$$ \rho_{2\rightarrow 1}\left( \eta_{12}  \left( \rho_{2\rightarrow 1}\left( \eta_{12}
	\left( \bullet 1(a) \oplus \bullet 2(b) \right) \right)
	\oplus \bullet 2(c) \right) \right)
$$
In what follows, MSOL($\tau_1$) stands for the monadic second order logic with quantification over subsets
of vertices. We denote by $G(\tau_1)$ the logic structure $<V(G),R>$ where $R$ is a binary
relation such that $R(u,v)$ holds if and only if $uv$ is an edge in $G$. An optimization problem is a
{\it LinEMSOL($\tau$) optimization problem} if it can be described in the following way (see \cite{LKLP}
for more details, since this is a version of the definition given by \cite{CMR} restricted to finite simple graphs),
$$  \opt \; \left\{   \sum_{1\le i \le l} a_i |X_i|   \;:\;
<G(\tau_1),X_1,\ldots,X_l > \; \vDash \theta(X_1,\ldots,X_l) \right\} $$
where $\theta$ is an MSOL($\tau_1$) formula that contains free set-variables
$X_1,\ldots,X_l, $ integers $a_i$ and \emph{Opt} is either $\min$ or $\max.$

We use the following result on LinEMSOL optimization problems.

\begin{thm} \emph{(Courcelle et all. \cite{CMR})} Let $k\in \mathbb{N}$ and let $\mathcal C$
be a class of graphs of clique-width at most $k$. Then every LinEMSOL$(\tau_1)$ optimization
problem on $\mathcal C$ can be solved in linear time if a $k$-expression of the graph is part
of the input.
\end{thm}

\noindent and we extend a result proved by Liedloff et al. (see Th. 31 in \cite{LKLP}) regarding
the complexity of the Roman domination decision problem to the corresponding decision problem
for the total Roman $\{2\}$-domination number.

\begin{thm} The total Roman $\{2\}$-domination problem is a LinEMSOL$(\tau_1)$ optimization problem.
\end{thm}

\noindent\textbf{Proof. } Let us show that the total Roman $\{2\}$-domination problem can be expressed
as a LinEMSOL$(\tau_1)$ optimization problem. Let $f=(V_0,V_1,V_2)$ be a total Roman $\{2\}$-domination
function in $G=(V,E)$ and let us define the free set-variables $X_i$ such that $X_i(v)=1$ whenever $v\in V_i$
and $X_i(v)=0,$ otherwise. For the sake of congruence with the logical system notation, we denote
by $|X_i|=\sum_{v\in V} X_i(v),$ even when, clearly, is $|X_i|=|V_i|.$
Observe that the total Roman $\{2\}$-domination decision problem corresponds to achieve the optimum
for the following expression.
$$  \min_{X_i} \; \left\{   |X_1|+2|X_2|   \;:\;
<G(\tau_1),X_0,X_1,X_2 > \; \vDash \theta(X_0,X_1,X_2) \right\} $$
where $\theta$ is the formula given by
$$ \begin{array}{l}
            \theta(X_0,X_1,X_2)= \left( \forall v \left( \left( X_1(v) \lor X_2(v) \right) \rightarrow \exists u \left( \left(X_1(u)\lor X_2(u)
            \right) \land R(u,v) \right) \right) \right) \land
     \\[.5em]
            \left( \forall v\left( X_1(v) \lor X_2(v) \lor \exists u \left(R(u,v)\land X_2(u) \right)\right.\right.
            \lor \; \left.\left. \exists u,w \left(R(u,v)\land R(w,v) \land X_1(u) \land X_1(w) \right) \right) \right).
    \end{array}
$$

Clearly, $\theta$ is an MSOL$(\tau_1)$ formula that describes the total Roman $\{2\}$-domination problem. Namely,
the formula has two main clauses. The first one requires that every vertex $v$ with a positive label $1$ or $2$
must have a neighbor $u$ with a positive label. The latter implies that the induced subgraph by the set of vertices
$V_1 \cup V_2$ has no isolated vertices. The second clause of the formula assures that for any vertex $v$ of the
graph either the vertex itself has a positive label, or either it has a neighbor with a label $2$, or either it has two
different neighbors having label $1$ each. Hence, when the formula $\theta$ is satisfied, the requirements of a total 
Roman $\{2\}$-domination assignment in $G$ holds. $\Box$

As a consequence, we may derive the following corollary

\begin{cor} The decision problem associated to the total Roman $\{2\}$-domination problem can be solved in linear time
on any graph $G$ with clique-width bounded by a constant $k$, provided that either there exists a linear-time
algorithm to construct a $k$-expression of $G$, or a $k$-expression of $G$ is part of the input.
\end{cor}

Since any graph with bounded treewidth is also a bounded clique-width graph, and it is well-known that any tree graph
has treewidth equal to 1, then we can deduce that the total Roman $\{2\}$-domination problem can be solved in linear
time for the class of trees. Besides, there are several classes of graphs $G$ with bounded clique-width $cw(G)$ like,
for example, the cographs ($cw(G)\le 2$) and the distance hereditary graphs ($cw(G)\le 3$), for which it is also
possible to solve the total Roman $\{2\}$-domination problem in linear time.

}


\section{Graphs with small or large total Roman $\{2\}$-domination}

As mentioned in Section 1, {for all ntc graphs }$G$ of order $n,$ $2\leq
\gamma_{tR2}(G)\leq n.$ {In this section, we }characterize all ntc graphs $G$
such that $\gamma_{tR2}(G)\in\{2,3,n\}$.

\begin{pro}
\label{TR2=2}Let $G$ be an ntc graph. {For any graph $H$, we have
$\gamma_{tR2} (K_{2}\vee H)=2$. Conversely, if $\gamma_{tR2}(G)=2$, there is
some graph $H$ such that $G=K_{2}\vee H$.}
\end{pro}

\noindent\textbf{Proof. }If $G=K_{2}\vee H$, then clearly $\gamma_{tR2}(G)=2$.
Conversely, assume that $\gamma_{tR2}(G)=2$ and let $f=(V_{0},V_{1},V_{2})$ be
a $\gamma_{tR2}(G)$-function. By definition of TR2DF of $G$, we have
$\gamma_{tR2}(G)=|V_{1}|+2|V_{2}|=2$. Since $G[V_{1}\cup V_{2}]$ has no
isolated vertex, we deduce that $|V_{2}|=0$ and $|V_{1}|=2$. Now, let
$V_{1}=\{x,y\}$. Clearly, $xy\in E(G),$ because $G[V_{1}\cup V_{2}]$ has no
isolated vertex, and every vertex in $V(G)-\{x,y\}$ is adjacent to both $x$
and $y$. Thus $G\cong K_{2}\vee H$, and $H$ is any graph of order $n-2$.
$\ \Box$

\begin{pro}
\label{TR2=3}Let $G$ be an ntc graph of order $n\geq5$
. Then $\gamma_{tR2}(G)=3$ if and only if either {$G$ has exactly one vertex
of degree $n-1$} or $\Delta(G)\leq n-2$ and $G$ is obtained from two disjoint
graphs $G_{1}$ and $G_{2}$ such that $G_{1}\in\{P_{3},C_{3}\}$ and $G_{2}$ is
any graph of order $n-3$ by adding edges between vertices of $G_{1}$ and
$G_{2}$ {in order that} every vertex of $G_{2}$ has at least two
neighbors in $G_{1}.$
\end{pro}

\noindent\textbf{Proof. }If $\Delta(G)=n-1$ and $G$ has exactly one vertex $u$
of degree $n-1$, then the function $f$ defined on $V(G)$ by $f(u)=2$, $f(v)=1$
for some $v\in V(G)-\{u\}$ and $f(w)=0$ for all $w\in V-\{u,v\}$ is a TR2DF
and so $\gamma_{tR2}(G)\leq3$. Since $G$ has exactly one vertex of degree
$n-1$, we deduce from Proposition \ref{TR2=2} that $\gamma_{tR2}(G)\geq3$ {and
the equality follows}$.$

Now assume that $\Delta(G)\leq n-2$ and $G$ is obtained from two disjoint
graphs $G_{1}$ and $G_{2}$ such that $G_{1}\in\{P_{3},C_{3}\}$ and $G_{2}$ is
any graph of order $n-3$ by adding edges between vertices of $G_{1}$ and
$G_{2}$ {in order that} every vertex of $G_{2}$ has at least two
neighbors in $G_{1}.$ Then the $f$ defined on $V(G)$ by $f(u)=1$ for all
vertex $u\in V(G_{1})$ and $f(v)=0$ for all $v\in V(G_{2})$ is a TR2DF of $G.$
Hence $\gamma_{tR2}(G)\leq3$, and the equality follows as above from
Proposition \ref{TR2=2}.

Conversely, assume that $\gamma_{tR2}(G)=3$. Then $G$ has at most one vertex
of degree $n-1$, for otherwise $\gamma_{tR2}(G)=2$ (by Proposition
\ref{TR2=2}). Let $f=(V_{0},V_{1},V_{2})$ be a $\gamma_{tR2}(G)$-function.
Since $\gamma_{tR2}(G)=|V_{1}|+2|V_{2}|=3,$ then it must be either
$|V_{1}|=|V_{2}|=1$ or either $|V_{1}|=3$ and $|V_{2}|=0.$ If $|V_{1}%
|=|V_{2}|=1$, with $V_{2}=\{u\}$ and $V_{1}=\{v\},$ then $uv\in E(G)$ and
every vertex in $V-\{x,y\}$ must be adjacent to $u$, because $f$ is a TR2DF.
So $u$ is the unique vertex of degree $n-1.$ Now assume that $V_{2}=\emptyset$
and $|V_{1}|=3$, where $V_{1}=\{u,v,w\}.$ Since $G[V_{1}]$ has no isolated
vertices, $G[V_{1}]\in\{P_{3},C_{3}\}$. Moreover, every vertex in
$V_{0}=V-\{u,v,w\}$ must be adjacent to at least two vertices of $V_{1}.$
Clearly, if $G_{1}=G[V_{1}]$ and $G_{2}=G[V_{0}],$ then $G$ is an ntc graph as
described in the statement. $\ \Box$

\begin{thm}
Let $G$ be an ntc graph of order $n$. Then $\gamma_{tR2}(G)=n$ if and only if
$G\in\{K_{2},K_{1,2}\}$ or every vertex of $G$ is either a leaf or a weak
support vertex.
\end{thm}

\textbf{Proof. }Assume that $\gamma_{tR2}(G)=n.$ Clearly, if $n\in\{2,3\},$
then $G\in\{K_{2},K_{1,2}\}.$ Hence assume that $n\geq4.$ Suppose first that
$G$ has a vertex $w$ which is neither a leaf nor a support vertex. Define the
function $f$ by $f(w)=0$ and $f(x)=1$ for all $x\in V(G)-\{w\}.$ Clearly $f$
is a TR2DF on $G$ with weight $n-1$, a contradiction. Thus, each vertex of $G$
is either a leaf or a support vertex. Now suppose that $G$ has a strong
support vertex, say $u.$ Let $u_{1}$ and $u_{2}$ be two leaves adjacent to
$u$. Define the function $f$ by $f(u)=2$, $f(u_{1})=f(u_{2})=0$ and $f(x)=1$
otherwise. Since $n\geq4,$ $f$ is clearly a TR2DF on $G$ of weight $n-1$, a
contradiction too. Therefore, every vertex of $G$ is either a leaf or a weak
support vertex as desired.

The converse is obvious. $\ \Box$

\section{Bounds}

We present in this section some bounds on the total Roman $\{2\}$-domination
number of ntc graphs in terms of the order, maximum and minimum degrees.

\begin{prop}
\label{upp1}Let $G$ be an ntc graph of order $n$ with girth $g\geq6$ and
minimum degree $\delta\geq2.$ Then $\gamma_{tR2}(G)\leq n+2-(\Delta+\delta).$
\end{prop}

\textbf{Proof. } Let $u$ be a vertex of $G$ of maximum degree and let $v$ be
any neighbour of $u$. Define the function $f$ on $V(G)$ by
$f(u)=1,f(v)=1,f(w)=0$ for all $w\in N(\{u,v\})-\{u,v\}$ and $f(w)=1$
otherwise. Since since $\delta\geq2$ and $g\geq6,$ set $A=V(G)-N[\{u,v\}]$ is
non-{empty} and no vertex of $A$ has two neighbors $N[\{u,v\}].$
Hence $f$ is well defined and is a TR2DF of weight $2+n-(\Delta+\deg_{G}(v))$,
and thus
\[
\gamma_{tR2}(G)\leq2+n-(\Delta+\deg_{G}(v))\leq n+2-(\Delta+\delta)
\]
$\ \Box$

{The sharpness of the previous bound can be seen by considering the cycles
}$C_{6}$ and $C_{7}.$ {Moreover, to see that }the condition $\delta\geq2$ is
essential in the statement of Proposition \ref{upp1}, consider the star
$K_{1,n-1}$ with $n\geq3,$ where $\gamma_{tR2}(K_{1,n-1})=3>n+2-(\Delta
+\delta)=2.$

{ }

\begin{cor}
Let $G$ be an ntc graph of order $n$ with girth $g\geq6$ and minimum degree
$\delta\geq2$ such that $\gamma_{tR2}(G)= n+2-(\Delta+\delta).$ Then for every
vertex $u$ of maximum degree we have that $d(v)=\delta(G)$ for all $v\in
N(u).$
\end{cor}

\begin{prop}
\label{upp2}Let $G$ be an ntc graph. Then
\[
\gamma_{tR2}(G)\geq\left\lceil \frac{2n}{\Delta+1}\right\rceil .
\]
If $\gamma_{tR2}(G) = \frac{2n}{\Delta+1}$ then $V_{2}=\emptyset$ for all
$\gamma_{tR2}(G)${-function $f=(V_{0},V_{1},V_{2})$.}
\end{prop}

\textbf{Proof. } Let $f=(V_{0},V_{1},V_{2})$ be a $\gamma_{tR2}(G)$-function
and let us denote by $V_{02}=\{w\in V_{0}:N(w)\cap V_{2}\neq\emptyset\}$ and
by $V_{01}=V_{0}-V_{02}.$ Thus $V(G)=V_{01}\cup V_{02}\cup V_{1}\cup V_{2}.$
Since any vertex $v\in V_{2}$ must have {at least one neighbor} in
$V_{1}\cup V_{2},$ {we deduce that $v\in V_{2}$} $|N(v)\cap
V_{02}|\leq\Delta-1$ and thus $|V_{02}|\leq(\Delta-1)|V_{2}|.$ Analogously,
$2|V_{01}|\leq(\Delta-1)|V_{1}|,$ because each vertex in $V_{01}$ must have at
least two neighbors in $V_{1}.$ Hence
\[%
\begin{array}
[c]{rcl}%
n & = & |V_{01}|+|V_{02}|+|V_{1}|+|V_{2}|\\[0.75em]
& \leq & \frac{\Delta-1}{2}|V_{1}|+(\Delta-1)|V_{2}|+|V_{1}|+|V_{2}|\\[0.75em]
& = & \frac{\Delta+1}{2}|V_{1}|+\Delta|V_{2}|\\[0.75em]
& \leq & \frac{\Delta+1}{2}\left(  |V_{1}|+2|V_{2}|\right)  =\frac{\Delta
+1}{2}\gamma_{tR2}(G)
\end{array}
\]
which leads to the desired result. { If $\gamma_{tR2}(G)= \frac{2n}{\Delta+1}$
then all the previous inequalities become equalities and hence $|V_{2}|=0.$
$\ \Box$ }

The sharpness of the bound in Proposition \ref{upp2} can be shown for cycles.

\section{Total Roman $\{2\}$-domination of trees}

In this section, we present a lower and upper bounds on the total Roman
$\{2\}$-domination number of trees. We start with a simple observation.

\begin{ob}
\label{strong} \emph{Let $G$ be a graph without isolated vertices and $v\in
V(G)$ a support vertex of $G$.  }

\begin{itemize}
\item \emph{For any total Roman $\{2\}$-dominating function $f$ of $G$,
$f(v)\ge1$.  }

\item \emph{If $v$ is a strong support vertex, then there exists a
$\gamma_{tR2}(G)$-function $f$ such that $f(v)=2$.  }
\end{itemize}

\end{ob}

\begin{thm}
\label{n+2-L}\emph{Let $T$ be a tree of order $n\geq2$ with $\ell(T)$ leaves.
Then
\[
\gamma_{tR2}(T)\geq\lceil\frac{2(n-\ell(T)+3)}{3}\rceil.
\]
} This bound is {sharp for paths}, stars and double stars.
\end{thm}

\textbf{Proof.} The proof is by induction on $n$. Clearly for all nontrivial
trees of order $n\leq4$ we have $\gamma_{tR2}(T)>\lceil\frac{2(n-\ell
(T)+3)}{3}\rceil$. For the inductive hypothesis, let $n\geq5$ and assume that
for every tree of order at least 2 and less than $n$ the result is true. Let
$T$ be a tree of order $n$. If $\mathrm{diam}(T)=2$, then $T$ is a star, which
yields $\gamma_{tR2}(T)=3=\lceil\frac{2(n-n+1+3)}{3}\rceil$. If $\mathrm{diam}%
(T)=3$, then $T$ is a double star and we have $\gamma_{tR2}(T)=4=\lceil
\frac{2(n-n+2+3)}{3}\rceil$. Henceforth we can assume $\mathrm{diam}(T)\geq4.$
Let $f$ be a $\gamma_{tR2}(T)$-function.

If $T$ has a strong support vertex $u$ with at least two leaves, say $u_{1}$
and $u_{2},$ then let $T^{\prime}=T-u_{1}$. By Observation \ref{strong},
$f(u)\geq1$ and we may assume without loss of generality that $f(u_{2})\geq
f(u_{1})$. Now the function $f$, restricted to $T^{\prime}$ is a TR2DF of
$T^{\prime}$ and we deduce from the inductive hypothesis that
\[
\gamma_{tR2}(T)=\omega(f)\geq\gamma_{tR2}(T^{\prime})\geq\lceil\frac
{2((n-1)-(\ell(T)-1)+3)}{3}\rceil=\lceil\frac{2(n-\ell(T)+3)}{3}\rceil.
\]
Thus in the sequel, we can assume that $T$ has no strong support vertex. Let
$v_{1}v_{2}\ldots v_{k}$ be a diametral path in $T$ and root $T$ in $v_{k}$.
Since $T$ has no strong support vertex, {any child} of $v_{3}$ is
a leaf or a support vertex of degree 2. We consider the following
cases.\newline

\smallskip\noindent\textbf{Case 1.} $\deg_{T}(v_{3})\geq3$.\newline First
suppose $v_{3}$ is a support vertex. By Observation \ref{strong}, we may
assume $f(v_{2})=f(v_{3})=2$. Let $T^{\prime}=T-v_{1}$ and define
$h:V(T^{\prime})\rightarrow\{0,1,2\}$ by $h(v_{2})=1$ and $h(x)=f(x)$ for
$x\in V(T^{\prime})-\{v_{2}\}$. Clearly $h$ is a TR2DF of $T^{\prime}$. It
follows from the induction hypothesis that
\[%
\begin{array}
[c]{lll}%
\gamma_{tR2}(T) & = & \omega(f)\\
& = & \omega(h)+1\\
& \geq & \gamma_{tR2}(T^{\prime})+1\\
& \geq & \lceil\frac{2((n-1)-\ell(T)+3)}{3}\rceil+1\\
& > & \lceil\frac{2(n-\ell(T)+3)}{3}\rceil,
\end{array}
\]
as desired. Now suppose $v_{3}$ is not a support vertex. Assume $u_{2}$ is a
child of $v_{3}$ and $u_{1}$ is a leaf adjacent to $u_{2}$. Clearly
$f(u_{1})+f(u_{2})\geq2$ and $f(v_{1})+f(v_{2})\geq2$. Assume without loss of
generality that $f(v_{2})\geq f(u_{2})$. Let $T^{\prime}=T-\{u_{1},u_{2}\}$.
If $f(v_{3})\geq1$ or $\deg(v_{3})\geq4$, then clearly the function $f$
restricted to $T^{\prime}$ is a TR2DF of $T$ and we conclude from the
inductive hypothesis that
\[%
\begin{array}
[c]{lll}%
\gamma_{tR2}(T) & = & \omega(f)\\
& = & \omega(f|_{T^{\prime}})+2\\
& \geq & \gamma_{tR2}(T^{\prime})+2\\
& \geq & \lceil\frac{2((n-2)-(\ell(T)-1)+3)}{3}\rceil+2\\
& > & \lceil\frac{2(n-\ell(T)+3)}{3}\rceil,
\end{array}
\]
as desired. Hence assume that $f(v_{3})=0$ and $\deg(v_{3})=3$. Let
$T^{\prime}=T-\{u_{1},u_{2},v_{1}\}$. Then the function $g:V(T^{\prime
})\rightarrow\{0,1,2\}$ defined by $g(v_{3})=1$ and $g(x)=f(x)$ for $x\in
V(T^{\prime})-\{v_{3}\}$, is a TR2DF of $T^{\prime}$ of weight $\gamma
_{tR2}(T)-2$. By the inductive hypothesis we have
\[%
\begin{array}
[c]{lll}%
\gamma_{tR2}(T) & = & \omega(f)\\
& = & \omega(g)+2\\
& \geq & \gamma_{tR2}(T^{\prime})+2\\
& \geq & \lceil\frac{2((n-3)-(\ell(T)-1)+3)}{3}\rceil+2\\
& > & \lceil\frac{2(n-\ell(T)+3)}{3}\rceil.
\end{array}
\]

\smallskip\noindent\textbf{Case 2.} $\deg_{T}(v_{3})=2$.\newline As above we
have $f(v_{1})+f(v_{2})\ge2$. If $f(v_{3})\ge1$, then the function
$g:V(T-v_{1})\to\{0,1,2\}$ defined by $g(v_{2})=1$ and $g(x)=f(x)$ for $x\in
V(T^{\prime})-\{v_{2}\}$, is a TR2DF of $T-v_{1}$ of weight $\gamma
_{tR2}(T)-1$ and by the inductive hypothesis we obtain
\[%
\begin{array}
[c]{lll}%
\gamma_{tR2}(T) & = & \omega(f)\\
& = & \omega(g)+1\\
& \ge & \gamma_{tR2}(T-v_{1})+1\\
& \ge & \lceil\frac{2((n-1)-\ell(T)+3)}{3}\rceil+1\\
& > & \lceil\frac{2(n-\ell(T)+3)}{3}\rceil.
\end{array}
\]

Hence let $f(v_{3})=0.$ If $f(v_{1})+f(v_{2})\geq3$, then reassigning
$v_{1},v_{2},v_{3}$ the value $1$ provides a $\gamma_{tR2}(T)$-function
$f^{\prime}$ for which $f^{\prime}(v_{3})\geq1,$ and this situation was
considered above. Therefore, we can assume that $f(v_{1})+f(v_{2})=2$. More
precisely, $f(v_{1})=f(v_{2})=1.$ It follows that $f(v_{4})\geq1.$ Let
$T^{\prime}=T-\{v_{1},v_{2},v_{3}\}.$ Clearly $T^{\prime}$ is nontrivial since
$\mathrm{diam}(T)\geq4.$ Now $T^{\prime}$ has order $2$, then $T$ is a path
$P_{5}$ and $\gamma_{tR2}(P_{5})=4\geq\lceil\frac{2(n-\ell(T)+3)}{3}\rceil.$
Hence suppose that $T^{\prime}$ has order at least three. Note that
$\ell(T)-1\leq\ell(T^{\prime})\leq\ell(T).$ Also, the function $f$ restricted
to $T^{\prime}$ is a T2RDF of $T^{\prime}$ of weight $\omega(f)-2.$ We deduce
from the inductive hypothesis on $T^{\prime}$ that
\[%
\begin{array}
[c]{lll}%
\gamma_{tR2}(T) & = & \omega(f)\\
& = & \omega(f|_{T^{\prime}})+2\\
& \geq & \gamma_{tR2}(T^{\prime})+2\\
& \geq & \lceil\frac{2((n-3)-\ell(T)+3)}{3}\rceil+2\\
& \geq & \lceil\frac{2(n-\ell(T)+3)}{3}\rceil,
\end{array}
\]
which competes the proof. $\ \ \Box$

\begin{lemm}
\label{lemma} \emph{If $T$ is a tree obtained from a path $v_{1}v_{2}\ldots
v_{k}\;(k\ge4)$ by adding a pendant path $v_{k-1}w$, then $\gamma_{tR2}(T)<
\frac{2(k+3)}{3}$.}
\end{lemm}

\textbf{Proof. } If $k\equiv0\pmod 3$, then define the function $f$ by
$f(v_{3i+1})=f(v_{3i+2})=1$ for $0\leq i\leq\frac{k}{3}-2$, $f(v_{k-2})=1$,
$f(v_{k-1})=2$ and {$f(v)=0$ for any remaining vertex $v$}. If
$k\equiv1\pmod 3$, then define the function $f$ by $f(v_{3i+1})=f(v_{3i+2})=1$
for $0\leq i\leq\frac{k-1}{3}-1$, $f(v_{k-1})=2$ and {$f(v)=0$ for
any remaining vertex $v$}. If $k\equiv2\pmod 3$, then define the function $f$
by $f(v_{3i+1})=f(v_{3i+2})=1$ for $0\leq i\leq\frac{k-2}{3}-1$,
$f(v_{k-2})=1$, $f(v_{k-1})=2$ and {$f(v)=0$ for any remaining
vertex $v$}. Clearly $f$ is an TR2DF of weight smaller than $\frac{2(k+3)}{3}%
$. $\hfill\Box$

{ }

\begin{thm}
\emph{For every tree $T$ of order }$n(T)\geq4$ with $s(T)$ support
vertices,\emph{
\[
\gamma_{tR2}(T)\leq\frac{3n(T)+2s(T)}{4}%
\]
with equality if and only if $T$ is the corona of a tree.}
\end{thm}

\textbf{Proof. }If $T$ is the corona of a tree $T^{\prime}$, then
$\gamma_{tR2}(T)={n(T)}=\frac{3n(T)+2s(T)}{4}.$ To prove that if
$T$ is a tree of order~$n(T)\geq4$ with $s(T)$ support vertices, then
$\gamma_{tR2}(T)\leq\frac{3n(T)+2s(T)}{4}$ with equality only if $T$ is the
corona of a tree, we proceed by induction on the order~$n(T)$. If $n(T)=4,$
then $T$ is either a star $K_{1,3},$ where $\gamma_{tR2}(K_{1,3}%
)=3<\frac{3n(T)+2s(T)}{4}$ or a path $P_{4}$ where $\gamma_{tR2}%
(P_{4})=4=\frac{3n(T)+2s(T)}{4}$ and $P_{4}$ is the corona of the path
$P_{2}.$ Let $n(T)\geq5$ and assume that every $T^{\prime}$ of order
$n(T^{\prime})<n(T)$ with $s(T^{\prime})$ support vertices satisfies
$\gamma_{tR2}(T^{\prime})\leq\frac{3n(T^{\prime})+2s(T^{\prime})}{4}$ with
equality only if $T^{\prime}$ is the corona of a tree. Let $T$ be a tree of
order $n(T)$. If $T$ is a star, then $\gamma_{tR2}(T)=3<\frac{3n(T)+2s(T)}{4}%
$. Likewise, if $T$ is a double star, then $\gamma_{tR2}(T)=4<\frac
{3n(T)+2s(T)}{4}$ (since $n(T)\geq5$). Henceforth, we can assume that $T$ has
diameter at least $4$. {Denote by {\normalsize $T_{x}$
the subtree induced by a vertex }$x$ and its
descendants in the rooted tree }$T.$

If $T$ has a strong vertex $u$ with {at least three leaves,} then
let $T^{\prime}$ be the tree obtained from $T$ by removing a leaf neighbor $w$
of $u.$ {Let }$f$ {be a }$\gamma_{tR2}(T^{\prime}%
)$-function $f$ such that $f(u)=2$, $f(v)\geq1$ for some $v\in N_{T^{\prime}%
}(u).$ {Clearly,} $f$ can be extended to TR2D-function of $T$ by
assigning a 0 to $w${, and thus }$\gamma_{tR2}(T)\leq\gamma
_{tR2}(T^{\prime}).$ {Now using }the induction on $T^{\prime}$
{and the fact that }$n(T^{\prime})=n(T)-1$ and $s(T^{\prime}%
)=s(T)$, {we obtain the desired result.} Henceforth, we can assume
that every support vertex of $T$ is adjacent to at most two leaves.

Let $v_{1}v_{2}\ldots v_{k}$ be a diametral path in $T$ such that $\deg
_{T}(v_{2})$ is as large as possible and root $T$ at $v_{k}$. Clearly
$\deg_{T}(v_{2})\in\{2,3\}.$ We consider the following cases.

\smallskip\noindent\textbf{Case 1.} $\deg_{T}(v_{2})=3$.\newline We
distinguish {the following subcases}.\newline\textbf{Subcase 1.1.}
$\deg_{T}(v_{3})\geq3$.\newline If $v_{3}$ is a support vertex or $v_{3}$ has
a child with degree 3 other than $v_{2}$, then any $\gamma_{tR2}(T-T_{v_{2}}%
)$-function can be extended to a TR2D-function of $T$ by assigning 2 to
$v_{2}$ and 0 to the leaf neighbors of $v_{2}$ and so $\gamma_{tR2}%
(T)\leq\gamma_{tR2}(T-T_{v_{2}})+2$. Since $T-T_{v_{2}}$ is a tree of order at
least four, by induction on $T-T_{v_{2}}$ and using the facts $n(T-T_{v_{2}%
})=n(T)-3$ and $s(T-T_{v_{2}})=s(T)-1$, we obtain $\gamma_{tR2}(T)\leq
\gamma_{tR2}(T^{\prime})+2\leq\frac{3n(T^{\prime})+2s(T^{\prime})}{4}%
+2<\frac{3n(T)+2s(T)}{4}$. Hence assume that every child of $v_{3}$ except
$v_{2}$ is of degree 2. Let $w_{2}$ be a child of $v_{3}$ besides $v_{2}$ and
let $w_{1}$ be the leaf neighbor of $w_{2}$. Clearly any $\gamma
_{tR2}(T-T_{w_{2}})$-function can be extended to a TR2D-function of $T$ by
assigning 1 to $w_{1},w_{2}$ and so $\gamma_{tR2}(T)\leq\gamma_{tR2}%
(T-T_{w_{2}})+2$. Note that $T-T_{w_{2}}$ is a tree of order at least four
with $n(T-T_{w_{2}})=n(T)-2$ and $s(T-T_{w_{2}})=s(T)-1$. Using the induction
on $T-T_{w_{2}}$, we obtain $\gamma_{tR2}(T)\leq\gamma_{tR2}(T^{\prime
})+2<\frac{3n(T^{\prime})+2s(T^{\prime})}{4}+2<\frac{3n(T)+2s(T)}{4}$.

\noindent\textbf{Subcase 1.2.} $\deg_{T}(v_{3})=2$ and $\deg_{T}(v_{4})\geq
3$.\newline If $T-T_{v_{3}}=P_{3}$, then clearly $\gamma_{tR2}(T)=5<\frac
{(3n(T)+2s(T))}{4}$. Hence assume that $T-T_{v_{3}}$ has order at least four.
Clearly any $\gamma_{tR2}(T-T_{v_{3}})$-function can be extended to a
TR2D-function of $T$ by assigning a $2$ to $v_{2}$, a $1$ to $v_{1}$ and a 0
to the leaves of $v_{2}$. It follows from the induction hypothesis on
$T-T_{v_{3}}$ and the facts $n(T-T_{v_{3}})=n-4$ and $s(T-T_{v_{3}})=s(T)-1$
that
\begin{align*}
\gamma_{tR2}(T) &  \leq\gamma_{tR2}(T-T_{v_{2}})+3\\
&  \leq\frac{3n(T-T_{v_{3}})+2s(T-T_{v_{3}})}{4}+3\\
&  \leq\frac{3(n(T)-4)+2(s(T)-1)}{4}+3<\frac{(3n(T)+2s(T))}{4}.
\end{align*}
\newline{\noindent\textbf{Subcase 1.3.} $\deg_{T}(v_{3})=2$ and $\deg
_{T}(v_{4})=2$.\newline}First let $\deg_{T}(v_{5})\geq3$. Hence $T-T_{v_{4}}$
has order at least three. If $T-T_{v_{4}}=P_{3}$, then it is easy to see that
$\gamma_{tR2}(T)=6<\frac{3n(T)+2s(T)}{4}$. Thus let $T-T_{v_{4}}\neq P_{3}$.
Then any $\gamma_{tR2}(T-T_{v_{4}})$-function can be extended to a
TR2D-function of $T$ by assigning a 2 to $v_{2}$ and $v_{3},$ and a $0$ to
other vertices in $T_{v_{4}}.$ Using the induction hypothesis on $T-T_{v_{4}}$
and the facts $n(T-T_{v_{4}})=n(T)-5$ and $s(T-T_{v_{4}})=s(T)-1$ we obtain
\begin{align*}
\gamma_{tR2}(T) &  \leq\frac{3n(T-T_{v_{4}})+2s(T-T_{v_{4}})}{4}+4\\
&  <\frac{3n(T)+2s(T)}{4}.
\end{align*}
\newline Assume now that $\deg_{T}(v_{5})=2$. {If $\deg(v_{i})\leq2$ for each
$i\geq5$, then the result follows from Lemma \ref{lemma}. Hence let $t$ be the
smallest integer such that $\deg(v_{t})\geq3$ for some $t\geq6$. }Let
$T^{\prime}=T-T_{v_{t-1}}.$ Note that $T^{\prime}$ has order at least three.
Suppose that $n(T^{\prime})=3,$ that is $T^{\prime}=P_{3}$\textbf{. }Then any
$\gamma_{tR2}(T_{v_{t-1}})$-function as defined in Lemma \ref{lemma} can be
extended to a TR2DF of $T$ by assigning a $2$ to $v_{t}$ and a $0$ to other
vertices of $T^{\prime},$ and clearly we have $\gamma_{tR2}(T)<\frac
{3n(T)+2s(T)}{4}$. Suppose now that $n(T^{\prime})\geq4.$ {If $t\equiv
1\pmod 3$, then any $\gamma_{tR2}(T-T_{v_{t-1}})$-function can be extended to
a TR2D-function of $T$ by assigning a }${2}${ to $v_{2}$, a }${1}${ to $v_{3}%
$, $v_{3i+2},v_{3i+3}$ for $1\leq i\leq\frac{t-1}{3}-1$ and a }${0}${ to the
remaining vertices of $T_{v_{t-1}.}$ Using the induction on $T-T_{v_{t-1}}$
and the fact }$\frac{2(t-1)}{3}$ {can be rewritten
}$\frac{3(t-1)}{4}-\frac{t-1}{12},$ {we have
\begin{align*}
\gamma_{tR2}(T) &  \leq\gamma_{tR2}(T-T_{v_{t-1}})+\frac{2(t-1)}{3}+1\\
&  \leq\frac{3n(T-T_{v_{4}})+2s(T-T_{v_{4}})}{4}+\frac{3(t-1)}{4}-\frac
{t-1}{12}+1\\
&  =\frac{3(n(T)-t)+2(s(T)-1)}{4}+\frac{3(t-1)}{4}-\frac{t-1}{12}+1\\
&  <\frac{3n(T)+2s(T)}{4}.
\end{align*}
Assume now that $t\equiv2\pmod 3$. Then any $\gamma_{tR2}(T-T_{v_{t-1}}%
)$-function can be extended to a TR2D-function of $T$ by assigning a }${2}${
to $v_{2}$, a }${1}${ to $v_{3i},v_{3i+1}$ for $1\leq i\leq\frac{t-2}{3}$ and
a }${0}${ to the remaining vertices of $T${$_{v_{t-1}}.$ By }the induction
hypothesis on $T-T_{v_{t-1}}$ we obtain
\begin{align*}
\gamma_{tR2}(T) &  \leq\gamma_{tR2}(T-T_{v_{t-1}})+\frac{2(t-2)}{3}+2\\
&  \leq\frac{3n(T-T_{v_{4}})+2s(T-T_{v_{4}})}{4}+\frac{3(t-2)}{4}-\frac
{t-2}{12}+2\\
&  =\frac{3(n(T)-t)+2(s(T)-1)}{4}+\frac{3(t-2)}{4}-\frac{t-2}{12}+2\\
&  <\frac{3n(T)+2s(T)}{4}.
\end{align*}
Finally, assume that $t\equiv0\pmod 3$. Then any $\gamma_{tR2}(T-T_{v_{t-1}}%
)$-function can be extended to a TR2D-function of $T$ by assigning a 2 to
$v_{2}$, a 1 to $v_{3}$, $v_{3i+1},v_{3i+2}$ for $1\leq i\leq\frac{t}{3}-1$
and a 0 to the remaining vertices of $T_{v_{t-1}}.$ By the induction
hypothesis on $T-T_{v_{t-1}}$ we have }%

\begin{align}
\gamma_{tR2}(T)  &  \leq\gamma_{tR2}(T-T_{v_{t-1}})+\frac{2t}{3}+1\nonumber\\
&  \leq\frac{3(n(T)-t)+2(s(T)-1)}{4}+\frac{3t}{4}-\frac{t}{12}+1\\
&  =\frac{3n(T)+2s(T)}{4}+\frac{6-t}{12}\leq\frac{3n(T)+2s(T)}{4}.\nonumber
\end{align}


If further $\gamma_{tR2}(T)=\frac{3n(T)+2s(T)}{4},$ then we have equality
throughout the previous inequality chain. In particular, we have $t=6$ and
$\gamma_{tR2}(T-T_{v_{t-1}})=\frac{3(n(T)-t)+2(s(T)-1)}{4}$. It follows from
the induction on $T-T_{v_{t-1}}$ that $T-T_{v_{t-1}}$ is the corona of a tree
and $v_{6}$ is support vertex (since $\deg_{T}(v_{6})\geq3$). It follows that
for any $\gamma_{tR2}(T-T_{v_{t-1}})$-function $g$, $g(v_{6})\geq1$ and
clearly $g$ can be extended to a TR2D-function of $T$ by assigning a $2$ to
$v_{2}$, a 1 to $v_{3},v_{5}$ and a 0 to other vertices in $T_{v_{5}}$. By the
induction hypothesis we obtain $\gamma_{tR2}(T)\leq\gamma_{tR2}(T-T_{v_{t-1}%
})+4<\frac{3n(T)+2s(T)}{4}$.

\smallskip\noindent\textbf{Case 2.} $\deg_{T}(v_{2})=2$.\newline By the choice
of the diametral path, we deduce that every child of $v_{3}$ with depth one
has degree two. Consider the following subcases.\newline\noindent
\textbf{Subcase 2.1.} $\deg_{T}(v_{3})\geq3$.\newline Suppose first that
$v_{3}$\textbf{ }is a strong support vertex, and let $u,w$ be two leaves of
$v_{3}.$ Let $T^{\prime}=T-\{u,v_{1},v_{2}\}.$ Clearly $T^{\prime}$ is a tree
of order $n(T^{\prime})=n(T)-3\geq4$ with $s(T^{\prime})=s(T)-1$ support
vertices. Let $g$ be a $\gamma_{tR2}(T^{\prime})$-function. Then we extend $g$
to a TR2D-function of $T$ by assigning a $1$ to $v_{1},v_{2}$ and a $0$ to
$u.$ In addition if $g(v_{3})\neq2$, then we reassign $v_{3}$ and $w$ the
{values $2$ and $0$} instead of $1$ to both. Now using the
induction hypothesis on $T^{\prime},$ we get%

\begin{align*}
\gamma_{tR2}(T)  &  \leq\gamma_{tR2}(T^{\prime})+2\\
&  \leq\frac{3n(T^{\prime})+2s(T^{\prime})}{4}+2\\
&  =\frac{3(n(T)-3)+2(s(T)-1)}{4}+2<\frac{3n(T)+2s(T)}{4}.
\end{align*}
Now, suppose that $v_{3}$ is not support vertex. Recall that every child of
$v_{3}$ is a support vertex of degree two. Let $T^{\prime}=T-T_{v_{3}}.$
Clearly $T_{v_{3}}$ has order $2\deg_{T}(v_{3})-1$ and $T^{\prime}$ has order
$n(T^{\prime})\geq2$ (since $\mathrm{diam}(T)\geq4$). If $n(T^{\prime})=2,$
then $\gamma_{tR2}(T)=$ $2\deg_{T}(v_{3})<\frac{3n(T)+2s(T)}{4},$ and if
$n(T^{\prime})=3,$ then $\gamma_{tR2}(T)=$ $2\deg_{T}(v_{3})+1<\frac
{3n(T)+2s(T)}{4}.$ Hence we assume that $n(T^{\prime})\geq4,$ and thus by
induction on $T^{\prime},$ $\gamma_{tR2}(T^{\prime})\leq\frac{3n(T^{\prime
})+2s(T^{\prime})}{4}.$ Since any $\gamma_{tR2}(T^{\prime})$-function
{can be extended} to a TR2D-function of $T$ by assigning a $0$ to
$v_{3}$ and a $1$ to each of the remaining vertices of $T_{v_{3}},$
$\gamma_{tR2}(T)\leq\gamma_{tR2}(T^{\prime})+2(\deg_{T}(v_{3})-1).$ Using the
fact that $s(T^{\prime})\leq s(T)-\deg_{T}(v_{3})+2,$ we obtain%

\begin{align*}
\gamma_{tR2}(T)  &  \leq\gamma_{tR2}(T^{\prime})+2(\deg_{T}(v_{3})-1)\\
&  =\frac{3n(T^{\prime})+2s(T^{\prime})}{4}+2(\deg_{T}(v_{3})-1)\\
&  \leq\frac{3(n(T)-2\deg_{T}(v_{3})+1)+2(s(T)-\deg_{T}(v_{3})+2)}{4}%
+2(\deg_{T}(v_{3})-1)\\
&  <\frac{3n(T)+2s(T)}{4}.
\end{align*}
Next we can assume that $v_{3}$ is a support vertex with $\deg_{T}(v_{3})=3$.
Let $T^{\prime}=T-\{v_{1},{v_{2}}\}$. As above we can easily see that
\begin{align*}
\gamma_{tR2}(T)  &  \leq\gamma_{tR2}(T^{\prime})+2\\
&  \leq\frac{3n(T^{\prime})+2s(T^{\prime})}{4}+2\\
&  \leq\frac{3n(T)+2s(T)}{4}.
\end{align*}
If further $\gamma_{tR2}(T)=\frac{3n(T)+2s(T)}{4},$ then we have equality
throughout the previous inequality chain. In particular, $\gamma
_{tR2}(T-\{v_{1},v_{2}\})=\frac{3(n(T)-2)+2(s(T)-1)}{4}$. It follows from the
induction on $T-\{v_{1},v_{2}\}$ that $T-\{v_{1},v_{2}\}$ is the corona of
some tree, implying that $T$ is the corona of a tree. \newline

\textbf{Subcase 2.2.} $\deg_{T}(v_{3})=2$ and $\deg_{T}(v_{4})\geq3$. \newline
If $T^{\prime}=T-T_{v_{3}}=P_{3}$, then clearly $\gamma_{tR2}(T)=5<\frac
{3n(T)+2s(T)}{4}$. Hence assume that $T^{\prime}\neq P_{3}$. If $v_{4}$ is
support vertex or has a child with depth 1 and degree at least 3, then clearly
there exists a $\gamma_{tR2}(T^{\prime})$-function that assigns a non-zero
positive value to $v_{4}$ and such a $\gamma_{tR2}(T^{\prime})$-function can
be extended to a TR2D-function of $T$ by assigning a $1$ to $v_{1},v_{2}$ and
a $0$ to $v_{3}.$ It follows from the induction hypothesis on $T^{\prime}$
that
\begin{align*}
\gamma_{tR2}(T) &  \leq\gamma_{tR2}(T^{\prime})+2\\
&  \leq\frac{3n(T^{\prime})+2s(T^{\prime})}{4}+2\\
&  <\frac{3n(T)+2s(T)}{4}.
\end{align*}
Now let $v_{4}$ has child $w_{2}$ with depth $\mathbf{1}$ and degree two, and
let $w_{1}$ be the leaf neighbor of $w_{2}$. Let $T^{\prime}=T-\{w_{2}%
,w_{1}\}.$ Clearly, $\gamma_{tR2}(T)\leq\gamma_{tR2}(T^{\prime})+2$. {By the
inductive hypothesis on }$T^{\prime}$ {and since }$T^{\prime}$ is not a
corona, $\gamma_{tR2}(T^{\prime})<\frac{3n(T^{\prime})+2s(T^{\prime})}{4}.$
{Using the facts that }$n(T^{\prime})=n(T)-2$ and $s(T^{\prime})=s(T)-1$ we
obtain
\begin{align*}
\gamma_{tR2}(T) &  \leq\gamma_{tR2}(T^{\prime})+2\\
&  <\frac{3n(T^{\prime})+2s(T^{\prime})}{4}+2\\
&  \leq\frac{3n(T)+2s(T)}{4}.
\end{align*}
Henceforth we assume that any child of $v_{4}$ is of depth 2. Thus $T_{v_{4}}$
is a tree obtain from a star by subdividing every edge twice. Let $w_{1}%
^{i}w_{2}^{i}w_{3}^{i}v_{4}$ be paths in $T$ where $w_{3}^{i}$
{"is a child} of $v_{4}$ for each $i\in\{1,2,\ldots,t\}$ and
$w_{3}^{1}=v_{3}$. If $t\geq3$, then any $\gamma_{tR2}(T-T_{v_{4}})$-function
can be extended to a TR2D-function of $T$ by assigning 1 to $v_{1},v_{2}%
,v_{3},v_{4},w_{2}^{i},w_{1}^{i}$ for $i\geq2$. Now we deduce from the
induction hypothesis on $T^{\prime}$ and the facts $n(T^{\prime})=n(T)-3t-1$
and $s(T^{\prime})\leq s(T)-t+1$ that
\begin{align*}
\gamma_{tR2}(T) &  \leq\gamma_{tR2}(T^{\prime})+2t+2\\
&  \leq\frac{3(n(T)-3t-1)+2(s(T)-t+1)}{4}+2t+2\\
&  <\frac{3n(T)+2s(T)}{4}.
\end{align*}
Hence assume that $t=2.$ If $\deg(v_{5})\geq3$, then let $T^{\prime
}=T-T_{v_{4}}.$ Then $\gamma_{tR2}(T)\leq\gamma_{tR2}(T-T_{v_{4}})+6.$ By the
induction hypothesis on $T^{\prime}$ and the facts $n(T-T_{v_{4}})=n(T)-7$ and
$s(T-T_{v_{4}})=s(T)-2$ we obtain
\begin{align*}
\gamma_{tR2}(T) &  \leq\gamma_{tR2}(T^{\prime})+6\\
&  \leq\frac{3(n(T)-7)+2(s(T)-2)}{4}+6\\
&  <\frac{3n(T)+2s(T)}{4}.
\end{align*}
Thus let $\deg(v_{5})=2$ and let $T^{\prime}=T-T_{v_{5}}.$ Note that
$T^{\prime}$ has order $n(T^{\prime})\geq2.$ If $n(T^{\prime})\in\{2,3\},$
then one can check that $\gamma_{tR2}(T)<\frac{3n(T)+2s(T)}{4}.$ Hence we
assume that $n(T^{\prime})\geq4.$ Then $\gamma_{tR2}(T)\leq\gamma
_{tR2}(T^{\prime})+6$. It follows from the induction hypothesis on
$T-T_{v_{5}}$ and the facts $n(T^{\prime})=n(T)-8$ and $s(T^{\prime})\leq
s(T)-1$ that
\begin{align*}
\gamma_{tR2}(T) &  \leq\gamma_{tR2}(T^{\prime})+6\\
&  \leq\frac{3(n(T)-8)+2(s(T)-1)}{4}+6\\
&  <\frac{3n(T)+2s(T)}{4}.
\end{align*}
\newline\newline\textbf{Subcase 2.3.} $\deg_{T}(v_{3})=\deg(v_{4})=2$.\newline
First let $\deg_{T}(v_{5})\geq3$. If $T^{\prime}=T-T_{v_{4}}=P_{3}$, then it
is easy to see that $\gamma_{tR2}(T)=5<\frac{3n(T)+2s(T)}{4}$. Hence assume
that $T^{\prime}\neq P_{3}$. If $v_{5}$ is a support vertex, then $v_{5}$ is
assigned a non-zero positive value under any {$\gamma_{tR2}(T^{\prime})$-set}
and thus one can easily see that $\gamma_{tR2}(T)\leq\gamma_{tR2}(T^{\prime
})+3.$ Using the induction hypothesis on $T^{\prime}$ and the facts
$n(T^{\prime})=n(T)-4$ and $s(T^{\prime})=s(T)-1$ we obtain
\begin{align*}
\gamma_{tR2}(T) &  \leq\frac{3n(T^{\prime})+2s(T^{\prime})}{4}+3\\
&  <\frac{3n(T)+2s(T)}{4}.
\end{align*}
If $v_{5}$ has child $w$ with depth one, then since there is a $\gamma
_{tR2}(T-T_{w})$-function that assigns a non-zero positive value to $v_{5},$
such a $\gamma_{tR2}(T-T_{w})$-function can be extended to a TR2D-function of
$T$ by assigning a $2$ to $w$ and 0 to other vertices in $T_{w}$. {By the
inductive hypothesis on }$T-T_{w}$ {and since }$T-T_{w}$ is not a corona,
$\gamma_{tR2}(T-T_{w})<\frac{3n(T-T_{w})+2s(T-T_{w})}{4}.$ {Moreover, we have
}$n(T-T_{w})\leq n(T)-2$ and $s(T-T_{w})=n(T)-1,$ and thus {%
\begin{align*}
\gamma_{tR2}(T) &  \leq\gamma_{tR2}(T-T_{w})+2\\
&  <\frac{3n(T-T_{w})+2s(T-T_{w})}{4}+2\\
&  \leq\frac{3(n(T)-2)+2(s(T)-1)}{4}+2\\
&  =\frac{3n(T)+2s(T)}{4}.
\end{align*}
} Suppose now that $v_{5}$ has child $w$ with depth two. Let $w$
{have $t_{3}$ leaves, $t_{2}$ children
with depth one and degree at least three and} $t_{1}$ children with depth one
and degree two. Let $T^{\prime}=T-T_{w}.$ Then any $\gamma_{tR2}(T^{\prime}%
)$-function can be extended to a TR2D-function of $T$ by assigning a $2$ to
every child of $w$ with depth one, $1+t$ to $w$ and 0 to other vertices in
$T_{w}$, where $t=0$ if $t_{3}=0$ and $t=1$ if $t_{3}\geq1$. {Clearly by the
inductive hypothesis on }$T^{\prime}$ {and since }$T^{\prime}$ is not a
corona, $\gamma_{tR2}(T^{\prime})<\frac{3n(T^{\prime})+2s(T^{\prime})}{4}.$
{Moreover, we know that }$n(T^{\prime})\leq n(T)-3t_{2}-2t_{1}-t_{3}-1$ and
$s(T^{\prime})=s(T)-t_{1}-t_{2}-t.$ {Now}

\begin{itemize}
\item {Assume that} $t_{2}\neq0$ {or} $t_{3}\neq0$. Then we have
\begin{align*}
\gamma_{tR2}(T)  &  \leq\gamma_{tR2}(T^{\prime})+2t_{2}+2t_{1}+1+t\\
&  <\frac{3n(T^{\prime})+2s(T^{\prime})}{4}+2t_{2}+2t_{1}+1+t\\
&  \leq\frac{3(n(T)-3t_{2}-2t_{1}-t_{3}-1)+2(s(T)-t_{1}-t_{2}-t)}{4}%
+2t_{2}+2t_{1}+1+t\\
&  \leq\frac{3n(T)+2s(T)}{4}.
\end{align*}
{ }

\item {Assume that} $t_{2}=0$ and $t_{3}=0$. Thus $t_{3}\geq1.$
{Using the fact that }there is a $\gamma_{tR2}(T^{\prime})$-function that
assigns a non-zero positive value to $v_{5},$ clearly then such a
$\gamma_{tR2}(T^{\prime})$-function can be extended to a TR2D-function of $T$
by assigning a $0$ to $w$ and $1$ to the remaining vertices of $T_{w}.$ It
follows that {%
\begin{align*}
\gamma_{tR2}(T)  &  \leq\gamma_{tR2}(T^{\prime})+2t_{1}\\
&  <\frac{3n(T^{\prime})+2s(T^{\prime})}{4}+2t_{1}\\
&  \leq\frac{3(n(T)-2t_{1}-1)+2(s(T)-t_{1})}{4}+2t_{1}\\
&  <\frac{3n(T)+2s(T)}{4}.
\end{align*}
}
\end{itemize}

Assume that $v_{5}$ has child with depth three and let $w_{1}w_{2}w_{3}%
w_{4}v_{5}$ be a path in $T$ where $w_{4}$ {is a child} of
$v_{5}$ different from $v_{4}$. {Considering the above cases} and
subcases we may assume that $\deg(w_{i})=2$ for $i\in\{1,2,3,4\}$. Clearly
$T-T_{w_{4}}$ has a $\gamma_{tR2}(T-T_{w_{4}})$-function $f$ such that
$f(v_{5})\geq1$, and $f$ can be extended to a TR2D-function of $T$ by
assigning a 1 to $w_{3},w_{2},w_{1}$ and 0 to $v_{4}$. Using the induction
hypothesis on $T-T_{w_{4}}$ we obtain
\begin{align*}
\gamma_{tR2}(T)  &  \leq\gamma_{tR2}(T-T_{w_{4}})+3\\
&  <\frac{3n(T-T_{w_{4}})+2s(T-T_{w_{4}})}{4}+3\\
&  <\frac{3n(T)+2s(T)}{4}.
\end{align*}
{Finally, assume }that $\deg_{T}(v_{5})=2$.
Let $f$ be $\gamma_{tR2}(T-T_{v_{3}})$-function such that $f(v_{4})$ is as
large as possible. It is easy to see that $f(v_{4})\geq1$ and $f$ can be
extended to a TR2D-function of $T$ by assigning a 1 to $v_{1},v_{2}$ and a $0$
to $v_{3}$. Using the induction hypothesis on $T-T_{v_{3}}$ we obtain
\begin{align*}
\gamma_{tR2}(T)  &  \leq\gamma_{tR2}(T-T_{v_{3}})+2\\
&  \leq\frac{3n(T-T_{v_{3}})+2s(T)}{4}+2\\
&  <\frac{3n(T)+2s(T)}{4}.
\end{align*}
$\Box$

{We conclude this section with two open problems.

\noindent {\bf Problem 1.} Is the problem of deciding whether $\gamma_{tR2}(G) =3\gamma(G)$ for a given
graph $G$  NP-hard.

\noindent {\bf Problem 2.} Characterize all graphs $G$ such that $\gamma_{tR2}(G)=3\gamma(G)$. }

\section*{Acknowledgments}
The authors thank the referees for their helpful comments and suggestions to improve the exposition and readability of the paper. H. Abdollahzadeh Ahangar was supported by the Babol Noshirvani University of Technology under research Grant Number BNUT/385001/98.


\begin{thebibliography}{9}                                                                                                %
\bibitem {bc}Bollob\'{a}s, B., Cockayne, E. J.: Graph-theoretic parameters
concerning domination, independence, and irredundance. J. Graph Theory 3
(1979) 241--249.

\bibitem {chha}M. Chellali, T.W. Haynes, S.T. Hedetniemi and A. MacRae, Roman
$\{2\}$-domination. \textit{Discrete Appl. Math.} \textbf{204} (2016) 22--28.

\bibitem{CL}H. Chen and C. Lu, A note on Roman ${2}$-domination problem in graphs.
\textit{arXiv:1804-09338}.

\bibitem {CDHH}E.J. Cockayne, P.A. Dreyer, S.M. Hedetniemi and S.T.
Hedetniemi, Roman domination in graphs. \textit{Discrete Math.} \textbf{278
}(2004) 11--22.

\bibitem{CMR}B. Courcelle, J.A. Makowsky and U. Rotics, Linear time solvable optimization
problems on graphs of bounded clique-width, {\textit Theory Comput. Syst.} {\bf 33} (2000) 125–150.

\bibitem {hv1}F. Harary and T.W. Haynes, Double domination in graphs,
\textit{Ars Combin.} \textbf{55} (2000), 201--213.

\bibitem {KM}W. Klostermeyer and G. MacGillivray, \textit{Roman, Italian, and
2-domination}, Manuscript, 2016.

\bibitem{LKLP}M. Liedloff, T. Kloks, J. Liu and S.L. Peng, Efficient algorithms for Roman
domination on some classes of graphs, \textit{Discrete Appl. Math.} \textbf{156} (2008) 34003415.

\bibitem{Re97} C. S. ReVelle, Can you protect the Roman Empire?
    \textit{Johns Hopkins Magazine} \textbf{49} (1997), no. 2, 40.

\bibitem {ReVelle}C.S. ReVelle and K.E. Rosing, \newblock Defendens imperium
romanum: a classical problem in military strategy. \textit{Amer. Math.
Monthly} 107 (7) (2000) 585-594.

\bibitem {Stewart}I. Stewart, \newblock Defend the Roman Empire!
\newblock\textit{Sci. Amer.} 281(6) (1999) 136-139.
\end{thebibliography}
\end{document}